\documentclass[12pt]{amsart}

\usepackage{amssymb}
\usepackage{times}
\usepackage{enumitem}

\usepackage[
  backend=biber,
  style=numeric,
  sorting=nyt,
  doi=false,
  isbn=false,
  url=false,
  giveninits=true,
]{biblatex}

\DeclareNameAlias{default}{given-family}
\DeclareFieldFormat[article]{title}{\mkbibitalic{#1}}
\DeclareFieldFormat[article]{journaltitle}{#1}
\renewbibmacro{in:}{}
\AtEveryBibitem{\clearfield{month}}
\AtEveryBibitem{\clearfield{number}}
\AtEveryBibitem{\clearlist{language}}
\DeclareFieldFormat{pages}{#1}
\renewbibmacro*{journal+issuetitle}{%
  \usebibmacro{journal}%
  \printtext{\addspace}%
  \printtext{\textbf{\thefield{volume}}}%
  \printtext{\addspace}%
  \printtext{(\thefield{year})}%
  \newunit%
}

\AtEveryBibitem{%
  \iffieldundef{shortjournal}{}%
  {\clearfield{journaltitle}%
   \clearfield{journaltitleaddon}%
   \renewbibmacro*{journal}{\printfield{shortjournal}}}%
}
\DeclareFieldFormat[book]{title}{#1}
\DeclareBibliographyDriver{book}{%
  \printnames{author}%
  \setunit{\addcomma\addspace}%
  \printfield{title}%
  \setunit{\addcomma\addspace}%
  \printlist{publisher}%
  \setunit{\addcomma\addspace}%
  \printlist{location}%
  \setunit{\addcomma\addspace}%
  \printfield{year}%
  \finentry%
}

\usepackage[colorlinks=true,linkcolor=blue,citecolor=blue]{hyperref}

\theoremstyle{plain}
\newtheorem{Thm}{Theorem}
\newtheorem{Cor}{Corollary}
\newtheorem{Lem}{Lemma}
\newtheorem{Prop}{Proposition}
\newtheorem{Rem}{Remark}

\newcommand{\thistheoremname}{}
\newtheorem*{genericthm*}{\thistheoremname}
\newenvironment{namedthm*}[1]
  {\renewcommand{\thistheoremname}{#1}\begin{genericthm*}}
  {\end{genericthm*}}

\theoremstyle{definition}
\newtheorem{Def}{Definition}
\newtheorem{Ex}{Example}

\newcommand{\R}{\mathbb{R}}
\newcommand{\td}{\mathrm{d}}
\DeclareMathOperator{\grad}{grad}
\DeclareMathOperator{\Hess}{Hess}
\DeclareMathOperator{\codim}{codim}
\numberwithin{equation}{section}

\begin{document}

\title{Concircular Sections Along Submanifolds}

\author{Morgan Robson}
\address{Morgan Robson\\
  Department of Mathematics\\
          South East Technological University\\
          Waterford\\
          Co. Waterford\\
          Ireland.}
\email{morgan.robson@setu.ie}

\date{\today}

\maketitle
\begin{abstract}
    We introduce an equation which generalises that of a concircular vector field to sections of the ambient tangent bundle along submanifolds. Totally umbilic hypersurfaces and various overdetermined systems (such as Obata's and Tashiro's equations) are obtained as specific cases. Given an ambient manifold $M$ and a connected submanifold $P$ of dimension at least $2$, the solution space $W(P;M)$ is shown to have dimension at most $\dim M+1$ with all solutions being of `special type' whenever it exceeds $\codim P+1$. This bound is shown to be sharp. When the submanifold $P$ is simply connected we show the maximum dimension is attained exactly when $P$ has constant sectional curvature $\kappa$ relative to the ambient manifold $M$. In the case of maximum dimension we introduce a bilinear form $\mathbb{I}$ on $W(P;M)$ which is Riemannian, degenerate or Lorentzian, depending on the sign of $\kappa$ and show $P$ immerses into $(W(P;M),\mathbb{I})$'s unit quadric. Hence our main result shows that whenever $P$ is simply connected, has dimension of at least $2$, and has constant curvature $\kappa\neq0$ relative to its ambient manifold $M$, there is an isometric immersion (which we construct explicitly) of $P$ with its metric induced from $M$, into the $\dim M$ - dimensional space form of curvature $\kappa$. 
\end{abstract}
\section{Introduction and Statement of Results}
Let $(M,g)$ be a Riemannian manifold and denote by $D$ the Levi-Civita connection of $g$. Given an immersed submanifold $P\subseteq M$, let $TM|_P$ denote the ambient tangent bundle of $M$ over $P$, i.e. $\pi^{-1}(P)$ where $\pi:TM \to M$ is the canonical projection. We study the pairs $(\sigma,\lambda)\in \Gamma(TM|_P)\times C^\infty(P)$ which solve the equation
     \begin{equation}\label{eq: focal point condition on base}
 D_{x}\sigma=\lambda x, \quad \forall x\in TP.
\end{equation}
Since $TM|_P=TP\oplus NP$, the section $\sigma$ need not be normal or tangent to $P$, although these two cases recover some well known scenarios. The normal case is easily recognised.
\begin{Ex}\label{ex: principal directions from split proportional distribution}
    If $P$ is an oriented hypersurface in $M$ and $\sigma \in \Gamma(NP)$ its unit normal, then equation (\ref{eq: focal point condition on base}) holds if and only if $P$ is totally umbilic in $M$.
\end{Ex}The tangent case is related to \textit{concircular vector fields}, introduced in \cite{fialkow_conformal_1939} and \cite{yano_concircular_1940}.
\begin{Def}\label{def: concircular vector field}
 Let $\Sigma$ be a pseudo-Riemannian manifold with Levi-Civita connection $\nabla$. The vector field $V\!\in\! \Gamma(T\Sigma)$ is called \textit{concircular} if
\[
\nabla_xV=\mu x, \quad \forall x\in T\Sigma,
\]
for some $\mu \in C^\infty(\Sigma)$. If there exists $\mathcal{H}\in C^\infty(\Sigma)$ and constants $\rho,c\in \R$ such that
\begin{align*}
&V=\grad \mathcal{H}, &&\text{and} &&\mu=\rho\mathcal{H}+c,
\end{align*}
then $V$ is called \textit{special concircular} \cite{tashiro_complete_1965}.
\end{Def}
\begin{Ex}\label{ex: concircular}
    Let $\sigma\in \Gamma(TP)$ be tangent to $P$. Equation (\ref{eq: focal point condition on base}) holds iff for all $x\in TP$,
    \begin{align}\label{eq: concircular equation}
    &D_x^\top\sigma=\lambda x &&\text{and} &&(D_x\sigma)^\bot=0,
    \end{align}
     which follows by projection onto $TP$ and $NP$. In the sense of Definition \ref{def: concircular vector field}, the leftmost equation is the concircular condition for $P$ equipped with its tangential connection, $D^\top$. By taking $P=M$ the rightmost equation becomes vacuous and the concircular condition for $M$ is recovered.
\end{Ex}
\begin{Ex}
    When the ambient manifold $M$ itself admits a concircular vector field, its restriction to $P$ satisfies equation (\ref{eq: focal point condition on base}).
\end{Ex}

Equation (\ref{eq: focal point condition on base}) is thus an `ambient', i.e. $TM$-valued generalisation of the concircular condition. We make the following definition.
\begin{Def}\label{def: special concircular}
    Let $(M,g)$ be a Riemannian manifold with $D$ the Levi-Civita connection of $g$ and $P$ an immersed submanifold. A section $\sigma\in \Gamma(TM|_P)$ is called \textit{concircular} if it satisfies equation (\ref{eq: focal point condition on base}). Furthermore $\sigma$ is called \textit{special} if
    \begin{align}\label{eq: special def for sections}
&\sigma^\top=\grad \mathcal{H}, &&\text{and} &&\lambda=\rho\mathcal{H}+c,
\end{align}
for $\sigma^\top$ the projection of $\sigma$ to $TP$, $\mathcal{H}\in C^\infty(P)$ and $\rho,c\in\mathbb{R}$ constants.
\end{Def}
Another point of contact with existing literature is the equation studied by He et al. in \cite{he_warped_2015}, which is itself a generalised Tashiro/Obata type equation \cite{tashiro_complete_1965}\cite{obata_certain_1962}. In \cite{he_warped_2015} the solution space to the equation $\Hess \lambda =\lambda q$ was studied, for $q$ a quadratic form on $P$ and $\lambda\in C^\infty(P)$. Equation (\ref{eq: focal point condition on base}) is similar in form, when solving for special solutions:
\begin{Ex}
 Let $(\sigma,\lambda)$ satisfy equation (\ref{eq: special def for sections}) for some $\mathcal{H}\in C^\infty(P)$ and $\rho,c\in\mathbb{R}$ with $\rho\neq 0$. Projecting equation (\ref{eq: focal point condition on base}) onto $TP$ and $NP$, and performing a linear change of variable gives
\begin{align*}\label{eq : T}
         \tag{T} &\Hess \lambda=\rho(\lambda g+\mathrm{II}),\\ 
         \tag{N} & (D_z\grad\lambda)^\perp=-\rho(D_z \sigma^\perp)^\bot \quad \forall z\in TP,
         \label{eq: N}
     \end{align*}
     where $\mathrm{II}:(x,y)\mapsto g(\sigma^\perp,D_xy)$ is the second fundamental form of $P$ in the direction $\sigma^\perp$, the normal part of $\sigma$. Equation (\ref{eq : T}) resembles that in \cite{he_warped_2015}, albeit with an additional term $\mathrm{II}$.
\end{Ex}
Following ideas in \cite{he_warped_2015} (and references therein), we analyse the solution space
\[
W(P;M)=\{(\sigma,\lambda)\in \Gamma(TM|_P)\times C^\infty(P): D_{x}\sigma=\lambda x,  \forall x\in TP\}.
\]
 When $\dim(P)=1$, equation (\ref{eq: focal point condition on base}) reduces to a first order linear ODE system for the components of $\sigma$, which can be solved locally for any $\lambda\in C^\infty(P)$. Hence $W(P;M)$ is infinite dimensional. When $\dim(P)\geq 2$ we will see that $W(P;M)$ is much smaller.\\ \\
Our results are as follows. Let $n:=\dim(M)$ and assume $\dim(P)\geq 2$ throughout. We start by showing that when $P$ is connected, $\dim~W(P;M)\leq n+1$  (Proposition \ref{prop: eval is injective}) and then investigate the zeros of solutions to equation (\ref{eq: focal point condition on base}), i.e. points $q\in P$ at which $\sigma_q=\mathbf{0}$. We show that any non-trivial solution of equation (\ref{eq: focal point condition on base}) has isolated zeros, and that two solutions are a constant multiple of one another whenever they share a zero (Proposition \ref{Prop: zero set is one dim}). We then consider the set of all zeros,
\begin{align*}
   &\mathcal{Z}:=\{q\in P : \exists(\sigma,\lambda)\in W(P;M)\backslash\{(\mathbf{0},0)\} \text{ such that } \sigma_q=\mathbf{0}\},
\end{align*}
and construct a map taking $q\in \mathcal{Z}$ to an element of $W(P;M)$ vanishing at $q$;
\[
\begin{aligned}
\mathcal{J}:&~\mathcal{Z} \to W(P;M),\\
&q \longmapsto (\sigma,\lambda) \quad \text{s.t.} \quad (\sigma|_q,\lambda(q))=(\mathbf{0},1).
\end{aligned}
\]
For $U\subseteq \mathcal{Z}$, an immersed submanifold of $P$, $\mathcal{J}|_U$ is shown to be a smooth immersion (Proposition \ref{prop: properties of Z}) and that either $\mathcal{Z}$ is nowhere dense in $P$ or $\dim W(P;M) \geq \dim P+1$ (Corollary \ref{cor: z nowhere dense}). We then show that a dimension bound on $W(P;M)$ implies all solutions to be locally special (cf. Definition \ref{def: special concircular}). When $P=M$, this is known, although it is stated with stronger assumptions than only a dimension bound (i.e. that at least two, everywhere pointwise linearly independent solutions exist - see \cite{shandra_geodesic_2019} and references therein). We show:
\begin{Thm}\label{thm: locally special}
    Suppose $P$ is connected, $\dim P \geq 2$ and $\dim W(P;M)\geq\mbox{codim}(P)+2$. Then there exists a constant $\rho\in\mathbb{R}$ such that
    \[
    \grad \lambda=\rho \sigma^\top,
    \]
    for all $(\sigma,\lambda)\in W(P;M)$. In particular every element of $W(P;M)$ is locally special.
\end{Thm}
We show that the space $W(P;M)$ splits when $M$ is a Riemannian product and $P$ a slice over a connected factor (Proposition \ref{prop: dimension split}). This is used to show the dimension bound in Theorem \ref{thm: locally special} is sharp. In other words, non-special concircular sections exist along submanifolds $P$ for which $\dim W(P;M)<\mathrm{codim}(P)+2$, in every codimension (Example \ref{ex: sharp}). Theorem \ref{thm: locally special} is then used to characterise when $\dim W(P;M)$ is maximal in terms of the sectional curvature of $M$, denoted $\mbox{sec}_M$, restricted to planes in $P$.
\begin{Def}\label{def: const sec curv rel M}
Let $R$ denote the Riemann curvature tensor of $M$. An immersed submanifold $P\subseteq M$ is said to have constant sectional curvature relative to $M$ if for some $\kappa\in \R$
    \begin{align}
        R(x,y)w=\kappa\big(g(y,w)x-g(x,w)y\big), \quad \forall x,y\in TP, \forall w\in TM|_P. \label{eq: I}
    \end{align}
    Equivalently, for every $x,y\in TP$ and $\nu\in NP$, $\sec_M(x,y)=\kappa$ and $R(x,y)\nu=0$.
\end{Def}
\begin{Thm}\label{thm: max dim W(P;M)}
Let $P$ be connected, immersed and $\dim P \geq 2$. If $\dim W(P;M)=n+1$ then $P$ has constant curvature relative to $M$. The converse holds if $P$ is simply connected.
\end{Thm}
This is well known in the case $P=M$ \cite{tashiro_complete_1965}. In particular, the maximum dimension bound is attained by any simply connected $P$ of dimension at least two, whenever $M$ has constant curvature (Corollary \ref{cor: max dim examples}). The case $\dim W(P;M)=n+1$ is studied in further depth. We show that for any two $(\sigma,\lambda),(\psi,\mu)\in W(P;M)$, the function 
\[
q\longmapsto \lambda(q)\mu(q)+\kappa g(\sigma,\psi)\big|_q,
\]
is constant on $P$, with $\kappa\in\mathbb{R}$ the constant value of $\sec_M$ for planes in $P$ (cf. Theorem \ref{thm: max dim W(P;M)}). We therefore define the following bilinear form on $W(P;M)$, 
\begin{align}\label{eq: I bilin form}
    \begin{aligned}
        \mathbb{I}:&W(P;M)\times W(P;M) \to \mathbb{R}\\
        &~\big((\sigma,\lambda),(\psi,\mu)\big) \longmapsto \lambda\mu+\kappa g(\sigma,\psi),
    \end{aligned}
\end{align}
which is Riemannian, degenerate, or Lorentzian when $\kappa>0$, $\kappa=0$ or $\kappa<0$ respectively. We remark that only $\dim W(P;M) \geq \codim P+2$ is needed for $\mathbb{I}$ to be sensibly defined, however we leave this direction of enquiry for future studies. We establish in the maximal dimensional case $\mathcal{Z}=P$ and thus prove the following.
\begin{Thm}\label{thm: J(P) is in quadric}
Suppose $P$ is connected, $\dim P\geq 2$ and $\dim W(P;M)=n+1$. The map $\mathcal{J}$ is an immersion of $P$ into $W(P;M)$'s unit quadric,
     \[
 \mathcal{J}(P)\subseteq Q:= \{(\sigma,\lambda)\in W(P;M):\mathbb{I}\big((\sigma,\lambda),(\sigma,\lambda)\big)=1\}.
 \]
 Furthermore, $\mathbb{I}$ restricted to $Q$ pulls back under $\mathcal{J}$ as
  \[
\mathcal{J}^\ast(\mathbb{I}|_Q)=\kappa(g|_P).
 \]
\end{Thm}
Hence combining Theorems \ref{thm: max dim W(P;M)} and \ref{thm: J(P) is in quadric} we get our main result.
\begin{Thm}\label{thm: isometric immersion}
    Suppose $P$ is simply connected, has constant sectional curvature $\kappa\neq0$ relative to $M$ and $\dim P \geq 2$. Then $(P,g|_P)$ immerses isometrically into the $n$-dimensional space form with curvature $\kappa$. 
\end{Thm}
Finally we make a brief remark on the origin of this work, to be explored at a later date. A solution $(\sigma,\lambda)\in W(P;M)$ can be lifted to a submanifold of $M$'s tangent bundle;
\[
\mathcal{G}:=\{(q,\sigma_q):q\in P\}\subseteq TM.
\]
The double tangent bundle of $M$ is split by $M$'s Levi-Civita connection into horizontal and vertical subspaces, $TTM=HM\oplus VM$, with their respective isomorphisms to $TM$
\begin{align*}
    &&\td\pi|_{HM}:HM \longrightarrow TM, &&\mathcal{K}|_{VM}:VM\longrightarrow TM.
\end{align*}
Submanifolds of $TM$ whose tangent vectors $\mathbb{X}$ solve $\mathcal{K}(\mathbb{X})=(\lambda\circ\pi)\td\pi(\mathbb{X})$ are modelled locally by submanifolds such as $\mathcal{G}$, i.e. the graphs of solutions to equation (\ref{eq: focal point condition on base}). When $M=\mathbb{S}^k$ ($k\in\mathbb{N}$), these solutions describe focal points of (possibly non-integrable) line families in $\mathbb{R}^{k+1}$, these arising as submanifolds of $T\mathbb{S}^k$. See \cite{georgiou_isotropic_2026} for further details.
\subsection*{Notation and Conventions} Throughout the following conventions are maintained.
\begin{enumerate}[label=\arabic*.]
    \item All manifolds are assumed in the smooth category.
    \item $M$ denotes a Riemannian $n$-manifold with metric $g$, Levi-Civita connection $D$ and curvature tensor $R$.
    \item $\Sigma\subseteq M$ is called immersed if $\Sigma$ is itself a manifold whose inclusion $\iota:\Sigma\hookrightarrow M$ is a smooth immersion. For an immersed submanifold we identify $T\Sigma$ with $\iota_\ast(T\Sigma)$.
    \item $P$ denotes a connected, immersed submanifold of $M$ with $\dim P \geq 2$. 
    \item  We denote $\pi^{-1}(P)$ i.e. the ambient tangent bundle of $M$ over $P$, by $TM|_P$.
\end{enumerate}
\subsection*{Acknowledgements}
The author would like to thank Nikos Georgiou and Brendan Guilfoyle for their helpful discussions and comments during the preparation of this paper.
\section{Fundamental Properties of the Solution Space and Theorem \ref{thm: locally special}}
Our first lemma relates elements of $W(P;M)$ and the curvature of the ambient manifold $M$. We remark that equation (\ref{eq: involutivity condition over base}) to follow is well known when $P=M$, i.e. in the usual concircular case \cite{chen_results_2020}.
\begin{Lem}\label{lem: Riemann curvature restriction}
    If $(\sigma,\lambda)\in W(P;M)$ then for all $x,y\in TP$,
        \begin{equation}\label{eq: involutivity condition over base}
        R(x,y)\sigma=\td\lambda(x)y-\td\lambda(y)x.
        \end{equation}
    Furthermore, by taking the trace over $TP$ (which we denote by $\mbox{tr}_P$)  we may solve for $
    \td\lambda$:
    \begin{equation}\label{eq: d tilde lambda constraint}
    \td\lambda=\tau(\cdot,\sigma),
    \end{equation}
    where $\tau\in \Gamma(\otimes^2(T M|_P)^\ast)$ is given by
    \[
    \tau(\mathbb{X},\mathbb{Y})=\frac{\mbox{tr}_P
    \big(Z\mapsto R(\mathbb{X},Z)\mathbb{Y}\big)}{{\dim(P)-1}},
    \]
    i.e. $\tau=-\frac{1}{\dim(P)-1}\mathrm{Ric}_P$, where $\mathrm{Ric}_P$ is the partial Ricci tensor over $P$.
\end{Lem} 
\begin{proof}
    The first equation is an easy consequence of equation (\ref{eq: focal point condition on base}), indeed for all $x,y\in TP$
    \begin{align*}
        R(x,y)\sigma&=D_xD_y\sigma-D_yD_x\sigma-D_{[x,y]}\sigma\\
        &=D_x(\lambda y)-D_y(\lambda x)-\lambda[x,y]\\
        &=\td\lambda(x)y-\td\lambda(y)x.
    \end{align*}
     The second follows as suggested, if $(e_i)_{i=1}^{\dim(P)}$ is a local orthonormal frame for $TP$
    \begin{align*}
        \tau(x,\sigma)=\frac{\sum\limits_{k=1}^{\dim P}g(R(x,e_k)\sigma,e_k)}{\dim(P)-1}
                    =\frac{\sum\limits_{k=1}^{\dim P}g(\td\lambda(x)e_k-\td\lambda(e_k)x,e_k)}{\dim(P)-1}
                    =\td\lambda(x).
    \end{align*}
\end{proof}
Since $\lambda$ is determined by the derivatives of $\sigma$ we employ the standard technique of prolongation \cite{branson_prolongations_2006} to equation (\ref{eq: focal point condition on base}). We work on the following Whitney sum
\[
E:=TM|_P\oplus (P\times \R),
\]
 so that $TM|_P$ encodes the possible values of $\sigma$ and the trivial bundle $P\times \R$ encodes the possible values of $\lambda$. We identify the spaces $C^\infty(P)$ and $\Gamma(P\times \R)$, thus consider solutions of equation (\ref{eq: focal point condition on base}) as sections of the bundle $E$, i.e.
\[
W(P;M)\subseteq \Gamma(TM|_P)\times \Gamma(P\times \mathbb{R})\simeq \Gamma(E).
\]
Furthermore equations (\ref{eq: focal point condition on base}) and (\ref{eq: d tilde lambda constraint}) comprise a linear system for $(\sigma,\lambda)$, hence elements of $W(P;M)$ correspond to parallel sections of $E$ for an appropriate choice of connection.
\begin{Def}
    Define the map $\widehat{\nabla}$ by
    \begin{equation}\label{eq: def of connection}
        \begin{aligned}
            \widehat{\nabla}:&~TP\times\Gamma(E) \to \Gamma(E)\\
            &(x,(w,f)) \longmapsto \big(D_xw-fx,\td f(x)-\tau(x,w)\big),
        \end{aligned}
    \end{equation}
    and as usual denote $\widehat{\nabla}(x,(w,f))$ by $\widehat{\nabla}_x(w,f)$.
\end{Def}
That $\widehat{\nabla}$ is a connection follows immediately since $\widehat{\nabla}=\nabla^\oplus-A$ where $A$ is the map
    \[
    \begin{aligned}
            A:&~TP\times\Gamma(E) \to \Gamma(E)\\
            &(x,(w,f)) \longmapsto (fx,\tau(x,w)),
    \end{aligned}
    \]
    and $\nabla^\oplus$ is $E$'s product connection defined by $\nabla^\oplus_x(w,f)=(D_xw,\td f(x))$ \cite{tu_differential_2017}. Therefore:
    \begin{Lem}\label{lem: sol iff par section}
        $(\sigma,\lambda)\in W(P;M)$, i.e. solves equation (\ref{eq: focal point condition on base}), iff $(\sigma,\lambda)$ is $\widehat{\nabla}$-parallel.
    \end{Lem}
    \begin{proof}
        \begin{align*}
            D_x\sigma=\lambda x, \quad \forall x\in TP& \iff D_x\sigma=\lambda x \quad\text{and}\quad \td\lambda(x)=\tau(x,\sigma), \quad \forall x\in TP\\
            &\iff\big(D_x\sigma-\lambda x, \td\lambda(x)-\tau(x,\sigma)\big)=(\mathbf{0},0), \quad \forall x\in TP\\
            &\iff \widehat{\nabla}_x(\sigma,\lambda)=(\mathbf{0},0), \quad \forall x\in TP.
        \end{align*}
    \end{proof}
    Hence we get an upper bound on the dimension of $W(P;M)$:
    \begin{Prop}\label{prop: eval is injective}
        Suppose $P$ is connected. For any $q\in P$ the evaluation map
        \begin{equation}
            \begin{aligned}
                \mathrm{eval}_q:&~W(P;M) \to T_q M \times \R\\
                &(\sigma,\lambda)\longmapsto (\sigma_q,\lambda(q)),
            \end{aligned}
        \end{equation}
        is injective. Hence $\dim W(P;M) \leq n+1$.
\end{Prop}
    \begin{proof}
    Elements of $W(P;M)$ are exactly the $\widehat{\nabla}$-parallel sections of $E$. Moreover, parallel transport on connected manifolds is uniquely determined by initial conditions.
    \end{proof}
The next proposition asserts that solutions to equation (\ref{eq: focal point condition on base}) always have isolated zeros and furthermore, that the location of these zeros distinguishes between solutions. Specifically, for $q\in P$, at most one $(\sigma,\lambda)\in W(P;M)$ (modulo scaling) satisfies $\sigma_q=\mathbf{0}$.
\begin{Prop}\label{Prop: zero set is one dim}
    Given $(\sigma,\lambda)\in W(P;M)$, the zeros of $\sigma$ are either all isolated or $\sigma \equiv \mathbf{0} $. Furthermore, for all $q\in P$ the subspace
    \[
    \mathcal{O}_q:=\{(\sigma,\lambda)\in W(P;M) : \sigma_q=\mathbf{0}\},
    \]
    has at most dimension $1$.
\end{Prop}
\begin{proof}
    To prove that the zeros are isolated, suppose $q\in P$ is a zero, i.e. $\sigma_q=\mathbf{0}$. If it happens that $\lambda(q)=0$ also, then $(\sigma,\lambda)_q=(\mathbf{0},0)$ and so $(\sigma,\lambda)$ agrees with the zero section of $E$ at $q$. Therefore by Proposition \ref{prop: eval is injective}, $(\sigma,\lambda)\equiv(\mathbf{0},0)$. Otherwise we may suppose $\lambda(q)\neq0$. In this case consider the derivatives of the function $\|\sigma\|^2\in C^\infty(P)$. If $x\in\Gamma(TP)$, then a quick application of equation (\ref{eq: focal point condition on base}) derives
    \begin{equation}\label{eq: derivative of sigma}
        D_x\|\sigma\|^2=2\lambda g(x,\sigma).
    \end{equation}
    Now differentiating again and applying equation (\ref{eq: focal point condition on base}) to remove the $D_x\sigma$ terms gives 
    \begin{align*}
        D_xD_x\|\sigma\|^2=2(D_x\lambda)\cdot g(x,\sigma)+2\lambda g(D_xx,\sigma)+2(\lambda)^2g(x,x).
    \end{align*}
Evaluating these equations at $q$ and taking $x$ such that $x_q\neq 0$ one finds (since $\sigma_q=\mathbf{0}$):
\begin{align*}
    &\|\sigma_q\|^2=0, && (D_x\|\sigma\|^2)_q=0, && (D_xD_x\|\sigma\|^2)_q=2(\lambda(q))^2g(x_q,x_q). 
\end{align*}
However by our assumptions on the values of $\lambda$ and $x$ at $q$, we have $(D_xD_x\|\sigma\|^2)_q>0$. Hence, since $x$ was arbitrary (up to being non-zero at $q$), $\|\sigma\|^2>0$ on some punctured neighbourhood of $q$, by Taylor's theorem to second order, and the zero is isolated as claimed. Now to prove that $\mathcal{O}_q$ is at most one dimensional, note that  $\mathrm{eval}_q(\mathcal{O}_q)\subseteq \{\mathbf{0}\}\times\R$, and since $\mathrm{eval}_q$ is injective, $\dim \mathcal{O}_q\leq \dim\left(\{\mathbf{0}\}\times\R\right)=1$.
\end{proof}
\begin{Rem}
    Due to the linearity of equation (\ref{eq: focal point condition on base}), the set $\mathcal{O}_q$ is non-trivial exactly when there are linearly independent elements in $W(P;M)$ which have linearly dependent vector components at $q$. Indeed if $(\sigma_i,\lambda_i)^{r}_{i=1}$ is a basis of $W(P;M)$
    \begin{align*}
        \exists(\sigma,\lambda)\in W(P;M)\backslash\{(\mathbf{0},0)\}: \sigma|_q=\mathbf{0} \iff\exists(c_1,\ldots,c_r)\in \mathbb{R}^r\backslash\{\mathbf{0}\} : \sum_{i=1}^r c_i\sigma_i|_q=\mathbf{0}.
    \end{align*}
\end{Rem}
We hence consider the following subset of $P$,  
\[
\mathcal{Z}:=\{q\in P : \mathcal{O}_q \neq \{\mathbf{0}\}\},
\]
and on this set define the following map
\begin{align}
\begin{aligned}
    \mathcal{J}:~&\mathcal{Z} \to W(P;M)\\
    &q \longmapsto (\mathrm{eval}_q)^{-1}(\mathbf{0},1).
\end{aligned}
\end{align}
In other words, given $q\in \mathcal{Z}$, the map $\mathcal{J}$ returns the non-trivial element $(\sigma,\lambda)\in W(P;M)$ such that $(\sigma|_q,\lambda(q))=(\mathbf{0},1)\in T_qM\times \mathbb{R}$. The existence and uniqueness of the target of $\mathcal{J}(q)$ follows from the definition of $\mathcal{Z}$ and Proposition \ref{prop: eval is injective} respectively. 
\begin{Prop}\label{prop: properties of Z}
    The set $\mathcal{Z}$ is closed in $P$. Furthermore whenever there exists $U\subseteq\mathcal{Z}$, with $U$ an immersed submanifold of $P$,
    \begin{enumerate}[label=\alph*)]
        \item The map $\mathcal{J}|_U:U\to W(P;M)$ is smooth.
        \item After identifying $W(P;M)$ with its tangent space, $\td\mathcal{J}|_U:TU \to W(P;M)$, the derivative of $\mathcal{J}|_U$ satisfies
        \begin{align*}
    \mathrm{eval}_q\big(\td\mathcal{J}|_U(q,x)\big)=(-x,0),
        \end{align*}
    for all $(q,x)\in TU$. In particular $\mathcal{J}|_U$ is an immersion of $U$ into $W(P;M)$.
    \item Every tangent vector in $T_qU$ arises as the value of some solution at $q$. I.e.
    \[
    T_qU\subseteq \{\sigma_q:(\sigma,\lambda)\in W(P;M)\}.
    \]
    \end{enumerate}
\end{Prop}
\begin{proof}
    In the following let $(\sigma_i,\lambda_i)_{i=1}^r$ be a basis of $W(P;M)$. A non-trivial linear combination of $(\sigma_i)_{i=1}^r$ vanishes at $p$ iff $\mathcal{O}_p\neq\{\mathbf{0}\}$. Let $(E_i)_{i=1}^n$ be a local frame for $TM|_P$ defined on an open neighbourhood $V$. Then
    \[
    (P\backslash\mathcal{Z})\cap V=\{p\in V: \big[\widehat{\boldsymbol{\sigma}}_1|\widehat{\boldsymbol{\sigma}}_2|\ldots|\widehat{\boldsymbol{\sigma}}_r\big](p)\in M_{n\times r}(\mathbb{R})\text{ has rank $r$}\},
    \]
    where $\widehat{\boldsymbol{\sigma}}_i$ are the component vectors of each $\sigma_i$ with respect to $(E_i)_{i=1}^n$. By continuity of the determinant it follows that $(P\backslash\mathcal{Z})\cap V$ and hence $P\backslash\mathcal{Z}$ is open, thus $\mathcal{Z}$ is closed. Now let $U$ be a subset of $\mathcal{Z}$, viewed as an immersed submanifold of $P$. To show $\mathcal{J}|_U$ is smooth write
    \[
    \mathcal{J}|_U : q\mapsto \sum_{i=1}^rc_i(q)(\sigma_i,\lambda_i),
    \]
    for maps $c_i:U\to\mathbb{R}$. By definition of $\mathcal{J}|_U$, $\mathrm{eval}_q(\mathcal{J}|_U(q))=(\mathbf{0},1)$ for all $q\in U$, i.e.
    \begin{align}\label{eq: J equations 1}
    &\sum_{i=1}^rc_i(q){\sigma_i}|_q=\mathbf{0},  &&\text{and}
    &&\sum_{i=1}^rc_i(q){\lambda_i}(q)=1.
    \end{align}
    As $(\sigma_i,\lambda_i)_{i=1}^r$ are linearly independent in $W(P;M)$ and $\mathrm{eval}_q$ is injective, $(\sigma_i|_q,\lambda_i|_q)_{i=1}^r$ are linearly independent in $T_qM\times\mathbb{R}$. After choosing (locally near each $q$) a non-singular $r\times r$ sub-matrix of the above system, smoothness of the $c_i$'s follows from Cramer's rule applied to this sub-matrix, since $\sigma_i$ and $\lambda_i$ are smooth. Hence $\mathcal{J}|_U$ is smooth. Now to calculate the derivative, differentiate equations (\ref{eq: J equations 1}) in the direction $x\in T_qU$;
    \begin{align}\label{eq: J equations 2}
    &\sum_{i=1}^r\td c_i(x){\sigma_i}|_q+c_i(q)\lambda_i(q)x=\mathbf{0},  &&\text{and}
    &&\sum_{i=1}^r\td c_i(x){\lambda_i}(q)+c_i(q)\td\lambda_i(x)=0,
    \end{align}
    where on the leftmost equation we have used equation (\ref{eq: focal point condition on base}). It follows from equations (\ref{eq: J equations 1}) that equations (\ref{eq: J equations 2}) simplify to
        \begin{align}\label{eq: J equations 3}
    &\sum_{i=1}^r\td c_i(x){\sigma_i}|_q=-x,  &&\text{and}
    &&\sum_{i=1}^r\td c_i(x)\lambda_i(q)=0,
    \end{align}
    with the rightmost equation following since by equation (\ref{eq: d tilde lambda constraint})
    \[
\sum_{i=1}^rc_i(q)\td\lambda_i(x)=\sum_{i=1}^rc_i(q)\tau(x,\sigma_i|_q)=\tau\left(x,\sum_{i=1}^rc_i(q)\sigma_i|_q\right)=\tau(x,\mathbf{0})=0.
    \]
    Thus
    \[
    \td\mathcal{J}|_U(q,x)=\sum_{i=1}^{r}\td c_i(x)(\sigma_i,\lambda_i).
    \]
    Evaluating both sides at $q\in U$ and using equations (\ref{eq: J equations 3}) completes the proof of (b), from which (c) follows.
\end{proof}
\begin{Cor}\label{cor: z nowhere dense}
    Either $\mathcal{Z}$ is nowhere dense in $P$ or $\dim W(P;M)\geq \dim P +1$
\end{Cor}
\begin{proof}
    If $\mathcal{Z}$ is not nowhere dense, then since $\mathcal{Z}$ is closed, there exists a non-empty open subset $U\subset \mathcal{Z}$. Since $\mathcal{Z}$ has the subset topology, $U$ is an open submanifold of $P$. Proposition \ref{prop: properties of Z} then implies that for all $q\in U$, for all $x\in T_qP$ 
    \[
    (x,0),(\mathbf{0},1)\in \mathrm{eval}_q(W(P;M)).
    \]
    Since the above elements span $T_qP\times\mathbb{R}$ and $\mathrm{eval}_q$ is injective, the dimension claim follows.
\end{proof}
We now move towards proving Theorem \ref{thm: locally special}. First we give two lemmas:
\begin{Lem}\label{lem: sigma closed}
    For any $(\sigma,\lambda),(\sigma',\lambda')\in W(P;M)$ we have
    \begin{align}\label{eq: wedge of the lambdas}
        \td\lambda\wedge(\sigma'^\top)^\flat+\td\lambda'\wedge(\sigma^\top)^\flat=0,
    \end{align}
    where $(\sigma^\top)^\flat$ is the one form on $P$ defined by $ x \longmapsto g(x,\sigma)$. Furthermore $(\sigma^\top)^\flat$ is closed.
\end{Lem}
\begin{proof}
To show $(\sigma^\top)^\flat$ is closed, let $x,y\in TP$. Then
\begin{align*}
\td(\sigma^\top)^\flat(x,y)&=x((\sigma^\top)^\flat(y))-y((\sigma^\top)^\flat(x))-(\sigma^\top)^\flat\big([x,y]\big)\\
&=x(g(y,\sigma))-y(g(x,\sigma))-g([x,y],\sigma)\\
&=0,
\end{align*}
by using the metric compatibility of $D$ and equation (\ref{eq: focal point condition on base}). Equation (\ref{eq: wedge of the lambdas}) follows since
\[
g(R(x,y)\sigma,\sigma')= \td\lambda(x)g(y,\sigma')-\td\lambda(y)g(x,\sigma')= (\td\lambda\wedge (\sigma'^\top)^\flat)(x,y),
\]
by equation (\ref{eq: involutivity condition over base}). The result now follows by
the anti-symmetry of the Riemann tensor;
    \[
    g(R(x,y)\sigma,\sigma')=-g(R(x,y)\sigma',\sigma).
    \]
\end{proof}
The above lemma shows that it is the tangential part of $\sigma$ which controls the derivative of $\lambda$, hence to prove Theorem \ref{thm: locally special} consider the subspaces
    \begin{align}\label{eq: t_q and n_q}
        &t_q:=\{\sigma^\top_q:(\sigma,\lambda)\in W(P;M)\}, && n_q:=\{\sigma^\perp_q:(\sigma,\lambda)\in W(P;M)\}
    \end{align}
    which track at each point how much of $W(P;M)$ lies tangential or normal to $P$.
\begin{Lem}\label{lem: points at which t_q has dimmmore than 2}
    The set
    \[
    \mathcal{U}:=\{q\in P:\dim t_q \geq 2\},
    \]
    is open. Furthermore, there exists a locally constant function $\rho\in C^\infty(\mathcal{U})$ such that
    \[
    \td\lambda =\rho(\sigma^\top)^\flat,
    \]
    on $\mathcal{U}$, for all $(\sigma,\lambda)\in W(P;M)$.
\end{Lem}
\begin{proof}
    If $\mathcal{U}$ is empty then it is open. Hence suppose $q\in \mathcal{U}$ and that $(\sigma_1,\lambda_1),(\sigma_2,\lambda_2)$ are elements in $W(P;M)$ such that $\sigma_1^\top|_q,\sigma_2^\top|_q\in T_qP$ are linearly independent. By continuity of the determinant (and by arguing similarly as at the start of Proposition \ref{prop: properties of Z}'s proof) $\sigma_1^\top$ and $\sigma_2^\top$ remain linearly independent on some open neighbourhood of $q$. Hence $\mathcal{U}$ is open. Now we derive the claimed equation. By applying equation (\ref{eq: wedge of the lambdas}) to the pairs  $(\sigma_1,\lambda_1),(\sigma_1,\lambda_1)$ and $(\sigma_2,\lambda_2),(\sigma_2,\lambda_2)$ we find
    \begin{align*}
        &\td\lambda_1 \wedge (\sigma_1^\top)^\flat=0, &&\text{and}  &&\td\lambda_2 \wedge (\sigma_2^\top)^\flat=0.
    \end{align*}
    As shown above there is an open neighbourhood of $q$, say $U\subseteq \mathcal{U}$ on which $\sigma_1^\top$ and $\sigma_2^\top$ are linearly independent, and hence non-zero. Therefore (after shrinking $U$ if needed) there exist functions $\rho_1,\rho_2\in C^\infty(U)$ such that
    \begin{align*}
        &\td\lambda_1|_U=\rho_1(\sigma_1^\top|_U)^\flat, &&\text{and}  &&\td\lambda_2|_U=\rho_2(\sigma_2^\top|_U)^\flat.
    \end{align*}
    If one now applies equation (\ref{eq: wedge of the lambdas}) to the pair $(\sigma_1,\lambda_1),(\sigma_2,\lambda_2)$, and uses the just derived expressions, we see
    \[
    (\rho_1-\rho_2)(\sigma_1^\top|_U)^\flat\wedge(\sigma_2^\top|_U)^\flat=0.
    \]
    Since $\sigma_1^\top|_U$ and $\sigma_2^\top|_U$ are linearly independent, so are $(\sigma_1^\top|_U)^\flat$ and $(\sigma_2^\top|_U)^\flat$, thus $\rho_1=\rho_2$. Denote this common function as simply $\rho$. It follows from $\td\lambda_i|_U=\rho(\sigma_i^\top|_U)^\flat$ that after taking an exterior derivative (and using that $(\sigma_i^\top|_U)^\flat$ is closed)
    \begin{align*}
        &0=\td \rho\wedge (\sigma_1^\top|_U)^\flat, &&\text{and}  &&0=\td \rho\wedge (\sigma_2^\top|_U)^\flat.
    \end{align*}
    which implies $\td\rho=0$ as $\sigma^\top_i|_U$ are linearly independent, i.e. $\rho$ is locally constant. Now we show that an arbitrary $(\sigma,\lambda)\in W(P;M)$ has the claimed form on $U$. On $U$ we note the following:
    \begin{align*}
    (\td\lambda-\rho(\sigma^\top)^\flat)\wedge(\sigma^\top_i)^\flat&=\td\lambda \wedge (\sigma_i^\top)^\flat-(\sigma^\top)^\flat\wedge \rho(\sigma_i^\top)^\flat\\
    &=\td\lambda \wedge (\sigma_i^\top)^\flat+\td\lambda_i\wedge(\sigma^\top)^\flat &&\text{(by $\td\lambda_i=\rho(\sigma_i^\top)^\flat$)}\\
    &=0, &&\text{(by equation (\ref{eq: wedge of the lambdas}))}
    \end{align*}
    for $i=1,2$. The claimed representation follows as the $\sigma_i^\top$'s are linearly independent on $U$. Since our chosen $q$ was arbitrary, we have the claimed representation locally on all of $\mathcal{U}$, with possibly different functions $\rho$ on each open neighbourhood $U$. However, it follows by $\td\lambda|_U=\rho(\sigma^\top|_U)^\flat$ on each $U$ that each $\rho$ agrees on overlaps, hence the function $\rho$ may be defined on all of $\mathcal{U}$, by local gluing.
\end{proof}
\begin{proof}[Proof of Theorem \ref{thm: locally special}]
    First we establish a dimension bound on $t_q$, for any $q\in P$. Consider the map
    \[
     \begin{aligned}
         T_q:&~W(P;M) \to T_qM\\
         &(\sigma,\lambda) \longmapsto \sigma_q.
     \end{aligned}
    \]
    The kernel of $T_q$ is exactly $\mathcal{O}_q$, and its range is contained in $t_q\oplus n_q$, these subspaces given by equation (\ref{eq: t_q and n_q}). Hence by rank nullity
    \begin{equation}\label{eq: dim formula}
        \dim W(P;M) \leq \dim(\mathcal{O}_q)+\dim(t_q)+\dim(n_q).
    \end{equation}
    From our dimension assumption, and since $n_q\subseteq N_qP$, and by Proposition \ref{Prop: zero set is one dim};
    \begin{align*}
    \dim(t_q)&\geq\dim W(P;M)-\dim(n_q)-\dim(\mathcal{O}_q)\\
    &\geq \mathrm{codim}(P)+2-\mathrm{codim}(P)-1\\
    &=1.
    \end{align*}
    Therefore $\dim (t_q) \geq 1$ for all $q\in P$. We now show that $\mathcal{U}$, as defined in Lemma \ref{lem: points at which t_q has dimmmore than 2} is dense in $P$, and then using the connectivity of $P$, that $\rho$ is constant on $\mathcal{U}$. Finally we show each $(\sigma,\lambda)$ has the claimed form.\\
    \textbf{1. Density.}\\
     Assume for a contradiction that $\mathcal{U}$ is not dense, i.e. there exists an open subset $W\subseteq P\backslash \mathcal{U}$. Then for all $q\in W$, $\dim(t_q)<2$, by definition of $W$. By equation (\ref{eq: dim formula})
    \[
    \dim(\mathcal{O}_q)\geq \dim W(P;M)-\dim(t_q)-\dim( n_q)>0,
    \]
Therefore $\mathcal{O}_q\neq\{\mathbf{0}\}$ for all $q\in W$, hence $W$ is an open subset of $\mathcal{Z}$. By Proposition \ref{prop: properties of Z} part (c), $t_q$ contains every tangent vector in $T_qW$, however $T_qW=T_qP$ as $W$ is open, which implies $\dim(t_q)=\dim(P)\geq 2$, contradicting the inequality for $\dim(t_q)$.\\
    \textbf{2. $\rho$ is constant on $\mathcal{U}$.}\\
    Now fix a basis for $W(P;M)$, say $\{(\sigma_i,\lambda_i)\}_{i=1}^r$. Define the function $\widetilde{\rho}\in C^\infty(P)$ by
    \[
    \widetilde{\rho}=\frac{\sum_{i=1}^r \td\lambda_i(\sigma_i^\top)}{\sum_{i=1}^rg(\sigma_i^\top,\sigma_i)}.
    \]
    This function is smooth since $\sigma_i,\lambda_i$ and $g$ are smooth, furthermore the denominator is never zero on $P$, as for any $q\in P$, $\dim t_q\geq 1$ implies that at least one $(\sigma_i,\lambda_i)$ satisfy 
    \[
    g(\sigma^\top_i,\sigma_i)|_q=\|\sigma_i^\top\|^2|_q\neq 0.
    \]
    It also follows from Lemma \ref{lem: points at which t_q has dimmmore than 2} that for any $q\in \mathcal{U}$
    \[
    \widetilde{\rho}(q)=\frac{\sum_{i=1}^r \td\lambda_i(\sigma_i^\top)|_q}{\sum_{i=1}^rg(\sigma_i^\top,\sigma_i)|_q}=\frac{\sum_{i=1}^r \rho(q)g(\sigma_i^\top,\sigma_i)|_q}{\sum_{i=1}^rg(\sigma_i^\top,\sigma_i)|_q}=\rho(q).
    \]
    Hence $\widetilde{\rho}|_\mathcal{U}=\rho$. Since $\mathcal{U}$ is open, $\td(\widetilde{\rho}|_\mathcal{U})=\td\rho=0$, with the last equality following since $\rho$ is locally constant. Hence since $\td\widetilde{\rho}$ vanishes on $\mathcal{U}$, a dense subset of $P$, and $\widetilde{\rho}$ is smooth, $\td\widetilde{\rho}$ vanishes on $P$. Hence $\widetilde{\rho}$ is locally constant. But $P$ is connected, thus $\widetilde{\rho}$ is constant on $P$. Hence $\rho=\widetilde{\rho}|_\mathcal{U}$ is a constant.\\
    \textbf{3. Claimed Form}\\
    For any $q\in P$, and any $x\in T_qP$, let $X\in \Gamma(TP)$ be a vector field such that $X|_q=x$. Also let $(q_n)_{n=1}^\infty$ be a sequence in $\mathcal{U}$ converging to $q$. Then for any $(\sigma,\lambda)\in W(P;M)$
    \[
    \td\lambda(x)=\lim_{n\to\infty}\td\lambda(X_{q_n})=\lim_{n\to \infty} \rho g(X_{q_n},\sigma)=\rho g(x,\sigma),
    \]
    with the second equality following from Lemma \ref{lem: points at which t_q has dimmmore than 2}. Hence $\grad \lambda=\rho \sigma^\top$. To finish, we remark that since $(\sigma^\top)^\flat$ is closed by Lemma \ref{lem: sigma closed}, locally there exists a smooth function $\mathcal{H}$ such that $\sigma^\top=\grad \mathcal{H}$. Hence $\grad \lambda=\rho\grad \mathcal{H}$, i.e.
    \begin{align*}
        \sigma^\top=\grad \mathcal{H}, &&\text{and} &&\lambda=\rho \mathcal{H}+c,
    \end{align*}
    for some constant $c$, locally.
    \end{proof}
    The above bound is sharp in the sense that (locally) non-special solutions can exist when $\dim W(P;M)\leq\mathrm{codim}(P)+1$. This is most easily seen in the codimension $0$ case, i.e. when $P=M$, so the bound for a counterexample becomes $\dim W(M;M)=1$. Indeed non special concircular vector fields exist and are plentiful \cite{chen_results_2015}. This can be used to construct non-zero codimensional examples. We first prove the following.
    \begin{Prop}\label{prop: dimension split}
        Suppose $\Sigma$ and $N$ are Riemannian manifolds with $\Sigma$ connected and $\nabla^\Sigma$ the Levi-Civita connection of $\Sigma$. Equip $\Sigma \times N$ with the product metric and Levi-Civita connection $D$. Let  $p\in N$ and set $P=\Sigma\times\{p\}$. Then
        \[
        W(P;\Sigma\times N)\simeq W(\Sigma;\Sigma) \times T_pN.
        \]
    \end{Prop}
    \begin{proof}
        By the definition of the product connection, $D_XY=\nabla^\Sigma_XY$ and  $D_XU=0$, for $X,Y\in \Gamma(T\Sigma)$, $U\in \Gamma(TN)$ \cite{oneill_semi-riemannian_1983}. Since $TP\simeq T\Sigma$ we identify vectors in these two bundles. Furthermore
        \[
        T(\Sigma\times N)|_P=T(\Sigma\times N)|_{\Sigma\times\{p\}}\simeq T\Sigma\oplus(\Sigma\times T_pN),
        \]
        with the last bundle being a Whitney sum. Hence if we let $(U_i)_{i=1}^{\dim(N)}$ denote a basis of $T_pN$, each section $\sigma\in \Gamma(T(\Sigma\times N)|_P)$ takes the form $\sigma=X+\sum_i\alpha_i U_i$ for $X\in \Gamma(T\Sigma)$ and $\alpha_i\in C^\infty(\Sigma)$ for each $i=1,\ldots,\dim(N)$. Therefore for all $x\in TP$:
        \begin{align*}
            D_x\sigma=\lambda x &\iff D_x (X+\sum_i\alpha_i U_i)=\lambda x,\\
            &\iff \nabla_x^\Sigma X+\sum_i(D_x\alpha_i)U_i=\lambda x, &&\text{(by product connection properties)}\\ 
            &\iff \nabla^\Sigma_x X=\lambda x \quad \text{and}\quad \sum_i(D_x\alpha_i)U_i=0, &&\text{(by projecting onto $T\Sigma$ and $T_pN$)}\\
            &\iff \nabla^\Sigma_x X=\lambda x \quad \text{and} \quad D_x\alpha_i=0, \quad \forall i. &&\text{(as the $U_i$ form a basis of $T_pN$)}
        \end{align*}
        Hence $(\sigma,\lambda)\in W(P;\Sigma\times N)$ if and only if $\sigma$ is the sum of an element of $W(\Sigma;\Sigma)$ (with corresponding function $\lambda$) and an element of $T_pN$ (the functions $\alpha_i$ being constant on $\Sigma$ by connectedness and the above). Hence the result follows.
    \end{proof}
    \begin{Ex}\label{ex: sharp}
        Let $\Sigma$, $N$ and $P$ be as in the previous proposition. We may take $\Sigma$ such that $\dim(\Sigma)\geq 2$ and $(V,\mu)\in W(\Sigma;\Sigma)$ which is locally non-special (see \cite{chen_results_2015} for existence). It follows from Theorem \ref{thm: locally special} that $\dim W(\Sigma;\Sigma)=1$, otherwise $(V,\mu)$ would be locally special. Hence
        \[
        \dim W(P;\Sigma\times N)=1+\dim(N) <\dim(N)+2=\mathrm{codim}(P)+2.
        \]
        Furthermore the lift of $V\in\Gamma(T\Sigma)$ to $\Gamma(T(\Sigma\times N))$  restricted to $P$ satisfies $(V,\mu)\in W(P;\Sigma \times N)$ and is a locally non-special section.
    \end{Ex}
\section{The Maximal Dimensional Case and Theorems \ref{thm: max dim W(P;M)} and \ref{thm: J(P) is in quadric}}
We now characterise when $\dim W(P;M)=n+1$, i.e. is maximal. We first calculate the curvature of the bundle $E$.
\begin{Prop}\label{prop: E curvature}
        The curvature tensor of $E$, denoted $R^{\widehat{\nabla}}$, is given by
        \begin{equation}
            R^{\widehat{\nabla}}(x,y)(w,f)=\big(R(x,y)w+\tau(y,w)x-\tau(x,w)y, (D_y\tau)(x,w)-(D_x\tau)(y,w)\big),
        \end{equation}
        where $x,y\in TP$ and $(w,f)\in\Gamma(E)$.
    \end{Prop}
    \begin{proof}
        We use the standard formula
        \[
        R^{\widehat{\nabla}}(x,y)(w,f)=\widehat{\nabla}_x\widehat{\nabla}_y(w,f)-\widehat{\nabla}_y\widehat{\nabla}_x(w,f)-\widehat{\nabla}_{[x,y]}(w,f).
        \]
        As $R^{\widehat{\nabla}}$ is tensorial in its $\Gamma(E)$ slot, it suffices to calculate its action on a local frame of $\Gamma(E)$. Let $(e_i)_{i=1}^{n}$ be a frame of $TM|_P$ and $x,y$ vector fields on $P$ all defined locally on an open subset. The set $\{(e_i,0)_{i=1}^n, (\mathbf{0},1)\}$ forms a local frame for $E$, where $\mathbf{0}$ is the zero section of $TM|_P$ and $0$, $1$ denote constant functions. Applying equation (\ref{eq: def of connection}) twice gives
        \begin{align*}
        \widehat{\nabla}_x\widehat{\nabla}_y(\mathbf{0},1)=(-D_xy,\tau(x,y))   
        \end{align*}
        and
        \begin{align*}
        \widehat{\nabla}_x\widehat{\nabla}_y(e_i,0)=\big(D_xD_ye_i+\tau(y,e_i)x,-\td(\tau(y,e_i))(x)-\tau(x,D_ye_i)\big).
        \end{align*}
        Furthermore the bracketed terms are
        \begin{align*}
        &\widehat{\nabla}_{[x,y]}(\mathbf{0},1)=\big(-[x,y],0\big)
        &&\text{and}
        &&\widehat{\nabla}_{[x,y]}(e_i,0)=\big(D_{[x,y]}e_i,-\tau([x,y],e_i)\big).
        \end{align*}
        Hence upon anti-symmetrising in $x$ and $y$ and subtracting the bracketed term, one finds
        \[
        R^{\widehat{\nabla}}(x,y)(\mathbf{0},1)=\big(-D_xy+D_yx+[x,y],\tau(x,y)-\tau(y,x)\big)=(\mathbf{0},0),
        \]
        which follows as $\tau$ is symmetric in its arguments due to the symmetry of the Riemann tensor (cf. equation (\ref{eq: d tilde lambda constraint})). Similarly, $R^{\widehat{\nabla}}(x,y)(e_i,0)=(A,B)$ where 
        \[
        A=R(x,y)e_i+\tau(y,e_i)x-\tau(x,e_i)y,
        \]
        and
        \[
        B=\{\td(\tau(x,e_i))(y)+\tau(y,D_xe_i)\}-\{\td(\tau(y,e_i))(x)+\tau(x,D_ye_i)\}+\tau([x,y],e_i).
        \]
        By definition of the induced connection on $\otimes^2(T M|_P)^\ast$
        \begin{align*}
            \td(\tau(y,e_i))(x)=(D_x\tau)(y,e_i)+\tau(D_xy,e_i)+\tau(y,D_xe_i),
        \end{align*}
        implying $B=(D_y\tau)(x,e_i)-(D_x\tau)(y,e_i)$. Hence putting this all together:
        \begin{align*}
        R^{\widehat{\nabla}}(x,y)(\mathbf{0},1)&=(\mathbf{0},0),\\
        R^{\widehat{\nabla}}(x,y)(e_i,0)&=\big(R(x,y)e_i+\tau(y,e_i)x-\tau(x,e_i)y,(D_y\tau)(x,e_i)-(D_x\tau)(y,e_i)\big).
        \end{align*}
        The claim now follows by the tensoriality of the curvature tensor.
    \end{proof}
\begin{proof}[Proof of Theorem \ref{thm: max dim W(P;M)}]
    We first prove the forward direction. Recall from Definition \ref{def: const sec curv rel M} that we must show for all $q\in P$, $x,y\in T_qP$ and $w\in T_qM$ that
        \begin{align*} \tag{\ref{eq: I}}
        R(x,y)w=\kappa\big(g(y,w)x-g(x,w)y\big).
    \end{align*}
    By our assumption on $W(P;M)$ and by Proposition \ref{prop: eval is injective}, $\mathrm{eval}_q:W(P;M)\to T_qM\times\mathbb{R}$ is an isomorphism. Hence there exists $(\sigma_w,\lambda_w)\in W(P;M)$ such that $\sigma_w|_q=w$. Furthermore,
    \[
    \dim W(P;M)=n+1=\mathrm{codim}(P)+\dim(P)+1 \geq \mathrm{codim}(P)+2,
    \]
    and as $P$ is connected, Theorem \ref{thm: locally special} applies. Hence
    \begin{align*}
        \td\lambda_w=\rho g(\cdot,\sigma_w),
    \end{align*}
    for some constant $\rho\in\mathbb{R}$. It now follows from equation (\ref{eq: involutivity condition over base}) that
    \[
    R(x,y)\sigma_w|_q=\td\lambda_w(x)y-\td\lambda_w(y)x=\rho\big(g(x,\sigma_w|_q)y-g(y,\sigma_w|_q)x\big),
    \]
    i.e.
    \[
    R(x,y)w=-\rho\big(g(y,w)x-g(x,w)y\big),
    \]
    as required, so that $\kappa=-\rho$. For the reverse direction, by Lemma \ref{lem: sol iff par section}
         \begin{align*}
             &\dim W(P;M)=n+1 \iff \text{$E$ has a frame of parallel sections}\iff \text{$E$ is $\widehat{\nabla}$-flat},
         \end{align*}
         with the last implication bi-directional since we assume $P$ is simply connected. Hence it suffices to show that the curvature of $E$ vanishes, i.e. by Proposition \ref{prop: E curvature}, that conditions (\ref{eq: i}) and (\ref{eq: ii}) hold: 
         \begin{align}
            \tag{i} R(x,y)w=\tau(x,w)y-\tau(y,w)x,\label{eq: i}\\
            \tag{ii} (D_y\tau)(x,w)=(D_x\tau)(y,w),\label{eq: ii}
        \end{align}
        for all $x,y\in TP$ and $w\in TM|_P$.
We have assumed that $P$ has constant sectional curvature relative to $M$, so that equation (\ref{eq: I}) with $\kappa\in \R$ holds. By the definition of $\tau$ (cf. equation (\ref{eq: d tilde lambda constraint}))
        \begin{align*}
            \tau(x,w)=\frac{\sum\limits_{k=1}^{\dim P}g\big(R(x,e_k)w,e_k\big)}{\dim P-1}=\frac{\kappa \sum\limits_{k=1}^{\dim P}g\big(g(e_k,w)x-g(x,w)e_k,e_k\big)}{\dim P-1}=-\kappa g(x,w),
        \end{align*}
        for all $x\in TP$, $w\in TM|_P$. Thus condition (\ref{eq: i}) follows immediately from equation (\ref{eq: I}) and condition (\ref{eq: ii}) can be seen to hold after a short calculation, using that $\kappa$ is constant, that the second fundamental form of $P$ is symmetric, and that $g$ is $D$-parallel.
    \end{proof}
\begin{Cor}\label{cor: max dim examples}
    If $M$ has constant sectional curvature and $P$ is a simply connected immersed submanifold with $\dim(P)\geq 2$, then $\dim W(P;M)=n+1$.
\end{Cor}
As we have seen in the proof of Theorem \ref{thm: max dim W(P;M)}, the constant $-\rho$ (as in Theorem \ref{thm: locally special}) is the constant sectional curvature of planes in $P$, relative to $M$, which we denoted by $\kappa$. Going forward we assume $\kappa$ denotes such a constant on $P$. 
\begin{Lem}\label{lem: I const}
    Fix $(\sigma,\lambda),(\psi,\mu)\in W(P;M)$. Given that $\dim W(P;M)=n+1$ and $P$ is connected, the function
    \[
    q\longmapsto \lambda(q)\mu(q)+\kappa g(\sigma,\psi)\big|_q,
    \]
    is constant on $P$.
\end{Lem}
\begin{proof}
    This is a straightforward calculation using Theorem \ref{thm: locally special}, and equation (\ref{eq: focal point condition on base}) and metric compatibility of $D$. For any $x\in\Gamma(TP)$:
    \begin{align*}
        D_x\big(\lambda\mu+\kappa g(\sigma,\psi)\big)&=\mu\td\lambda(x)+\lambda\td\mu(x)+\kappa D_x(g(\sigma,\psi))\\
        &=-\kappa\mu g(x,\sigma)-\kappa\lambda g(x,\psi)+\kappa g(D_x\sigma,\psi)+\kappa g(\sigma,D_x\psi)\\
        &=-\kappa\mu g(x,\sigma)-\kappa\lambda g(x,\psi)+\kappa g(\lambda x,\psi)+\kappa g(\sigma,\mu x)\\
        &=0.
    \end{align*}
    Hence the function is locally constant. As $P$ is connected, it is constant.
\end{proof}
We are therefore able to define the bilinear form $\mathbb{I}$ on $W(P;M)$ given in equation (\ref{eq: I bilin form}), by defining $\mathbb{I}\big((\sigma,\lambda),(\psi,\mu)\big)$ to be equal to the constant value of the function in Lemma \ref{lem: I const}. We remark that $\mathbb{I}$ can sensibly be defined whenever $\dim W(P;M)\geq \codim P+2$, as an identical argument to that given in Lemma \ref{lem: I const} with `$-\rho$' in place of $`\kappa$' shows the consistency of the required function. We are now ready to prove Theorems \ref{thm: J(P) is in quadric} and \ref{thm: isometric immersion}.
\begin{proof}[Proof of Theorem \ref{thm: J(P) is in quadric}]
    First note that since for all $q\in P$, $\mathrm{eval}_q$ is bijective in the maximal dimension case, $(\mathrm{eval}_q)^{-1}(\mathbf{0},1)\in\mathcal{O}_q$. Hence $\mathcal{Z}=P$ and so $\mathcal{J}$ immerses $P$ by Proposition \ref{prop: properties of Z}. To see that $\mathcal{J}(P)$ is contained in the quadric, let $q\in P$ and simply note that by definition, $\mathcal{J}(q)=(\sigma,\lambda)$ for the unique $(\sigma,\lambda)$ which evaluates to $(\mathbf{0},1)$ at $q$. Therefore
    \[
    \mathbb{I}\big(\mathcal{J}(q),\mathcal{J}(q)\big)=\lambda^2+\kappa \|\sigma\|^2=1^2+\kappa\cdot 0=1,
    \]
    where the last equality follows by evaluating at $q$ and by Lemma \ref{lem: I const}. Now to show the claimed relationship between the metrics, note that by Proposition \ref{prop: properties of Z}, if $X,Y\in T_qP$, then $\td\mathcal{J}(q,X)=(\sigma,\lambda)$ and $\td\mathcal{J}(q,Y)=(\psi,\mu)$ for some elements of $W(P;M)$ such that $(\sigma|_q,\lambda(q))=(-X,0)$ and $(\psi|_q,\mu(q))=(-Y,0)$.
    Hence
    \[
    \mathbb{I}(\td\mathcal{J}(q,X),\td\mathcal{J}(q,Y))=\big(\lambda\mu+\kappa g(\sigma,\psi)\big)\Big|_q=\kappa g(X,Y).
    \]
\end{proof}
\begin{proof}[Proof of Theorem \ref{thm: isometric immersion}]
Since $P$ is simply connected and has constant curvature $\kappa$ relative to $M$, Theorem \ref{thm: max dim W(P;M)} implies that $\dim W(P;M)=n+1$. Hence Theorem \ref{thm: J(P) is in quadric} implies $P$ immerses into $Q$. Equip $W(P;M)$ with the rescaled metric $\widetilde{\mathbb{I}}~=\kappa^{-1}\mathbb{I}$ so that
\[
Q=\{(\sigma,\lambda)\in W(P;M): \widetilde{\mathbb{I}}\big((\sigma,\lambda),(\sigma,\lambda)\big)=\kappa^{-1}\}.
\]
Therefore $Q$ is the sphere of radius $\kappa^{-1/2}$ in $(W(P;M),\widetilde{\mathbb{I}})$ which is identified with either Euclidean space (if $\kappa>0$) or Minkowski space (if $\kappa<0$).
Any connected component of $Q$ is a model of the $n$-dimensional space form with curvature $\kappa$. As $P$ is connected, $\mathcal{J}(P)$ lies in a single connected component. The result follows by Theorem \ref{thm: J(P) is in quadric}.
\end{proof}
\subsection*{AI Disclosure}
No mathematical content or text in this paper was generated by LLM usage. LLM were used to help locate relevant literature and for grammatical/typographical proof reading. The author retains full responsibility for the accuracy of the work.
\printbibliography
\end{document}